\documentclass{amsart}
\usepackage{amsfonts,amssymb,amsthm}
\usepackage[utf8]{inputenc}
\usepackage{mathrsfs}
\usepackage{enumerate}
\usepackage{enumitem}
\usepackage{upref}
\usepackage[all]{xy}

\theoremstyle{definition}
\newtheorem{defin}{Definition}[section]
\theoremstyle{plain}
\newtheorem{thm}[defin]{Theorem}

\newtheorem{prop}[defin]{Proposition}
\newtheorem{lemma}[defin]{Lemma}
\newtheorem{step}{Step}
\theoremstyle{remark}
\newtheorem{ex}[defin]{Example}
\newtheorem{remark}[defin]{Remark}

\renewcommand{\qedsymbol}{\ensuremath{\blacksquare}}

\begin{document}

\title{Laudal's Lemma in positive characteristic}
\author{Paola Bonacini}
\address{
University of Catania\\
Viale A. Doria 6\\
95124, Catania, Italy
}
\email{bonacini@dmi.unict.it}

\begin{abstract}
Laudal's Lemma states that if $C$ is a curve of degree $d>s^2+1$ in $\mathbb{P}^3$ over an algebraically closed field of characteristic 0 such that its plane section is contained in an irreducible curve of degree $s$, then $C$ lies on a surface of degree $s$. We show that the same result does not hold in positive characteristic and we find different bounds $d>f(s)$ which ensure that $C$ is contained in a surface of degree $s$.
\end{abstract}

\maketitle
\pagestyle{plain}

\section{Introduction}

Let $C$ be a curve in $\mathbb{P}^3_k$, being $k$ an algebraically closed field. Let $\Gamma$ be the generic plane section of $C$. In this paper we study the problem of finding bounds on the degree of $C$ in such a way that, if $\Gamma$ is contained in a plane curve of degree $s$, then $C$ is contained in a surface of the same degree. In the case that char\,$k=0$ the following result has been proved:

\begin{thm}[Laudal's Lemma, {\cite[Corollary, p.147]{Lau}},\cite{GruPes}]  \label{T:2}
If $\Gamma$ is contained in an integral plane curve of degree $s$ and $\deg C>s^2+1$, then $C$ is contained in a surface of degree $s$.
\end{thm}

The bound on the degree of the curve found in this result is sharp. Indeed, there are examples of curve of degree $s^2+1$ whose the generic plane section is contained in an irreducible plane curve of degree $s$ and that are not contained in any surface of degree $s$ (see \cite{Har}, \cite{GruPes} and \cite[Proposition 1]{Strano}).

In this paper, following the proof of Gruson and Peskine of Laudal's Lemma in \cite{GruPes}, we prove an analogous result in the case that the field $k$ has positive characteristic:
\begin{thm}   \label{T:1}
Let $C \subset \mathbb P^3$ be a non degenerate reduced curve of degree $d$ in characteristic $p>0$. Suppose that $\Gamma$ is contained in an integral plane curve of degree $s$. Then $C$ is contained in a surface of degree $s$, if one of the following conditions is satisfied:
\begin{enumerate}
	\item $C$ is connected, $p\ge s$ and $d>s^2+1$;
	\item $C$ is connected, $p<s$ and $d>s^2+p^{2n}$, with $p^n<s\le p^{n+1}$; in particular this holds if $d>2s^2-2s+1$;
	\item $p>s$ and $d>s^2+1$;
	\item $p\le s$ and $d>s^2+p^{2n}$, with $p^n\le s<p^{n+1}$. In particular this holds if $d>2s^2$.
\end{enumerate}
\end{thm}

Let us make a note about terminology. Given the incidence variety $T=\{([H],P)\in \check{\mathbb P}^3\times \mathbb P^3\mid P\in C\cap H\}$ associated to $C$, the fibre of the projection $T\rightarrow \check{\mathbb P^3}$ over the generic point $\eta\in \check{\mathbb P^3}$ is the generic plane section $\Gamma$. In particular we consider the open subset $U\subset \check{\mathbb P}^3$ such that any $[H]\in U$ corresponds to a plane $H\subset \mathbb P^3$ that $C$ meets transversally and such planes are generic for the curve $C$.

Let us give now a sketch of the proof of Theorem \ref{T:1}, given in Section 4. We follow the idea of the proof of Theorem \ref{T:2} given by Gruson and Peskine in \cite{GruPes}. So we take $S\subset \check{\mathbb P}^3\times \mathbb P^3$ containing $T$ such that the fibre over $\eta$ is an integral plane curve of degree $s$ containing $\Gamma$. Then we suppose that $h^0(\mathscr I_C(s))=0$ and, using Theorem \ref{T:21}, that is the main result of Section 3, we factor the projection  $S\rightarrow \mathbb P^3$ through a generically smooth morphism $S_r\rightarrow \mathbb P^3$, with $S_r=S\times_{\mathbb P^3,\,F^r} \mathbb P^3$ and $F^r$ some $r$-th power of the absolute Frobenius of $\mathbb P^3$. Then, proceeding as in \cite{GruPes}, we arrive to the inequality $d\le s^2+p^{2r}$. Remarking that it must be $h^1(\mathscr I_C(s-p^r))\ne 0$ we find the desired inequalities.

Looking at the proof we see that the assumption that $C$ is reduced is required to find a suitable bound to the power $p^r$. Moreover, in the case that $C$ is connected this bound is sharp, as we will see in Example \ref{E}. Indeed, generalizing the example given in \cite{GruPes} and \cite[Proposition 1]{Strano} to prove that the bound in Theorem \ref{T:2} is sharp, we consider the sheaf $\mathscr E={F^n}^{\star}(\mathscr E_0)$, with
$\mathscr E_0$ null-correlation bundle and $F$ absolute Frobenius on $\mathbb P^3$. Then
the zero locus of a generic global section of $\mathscr E(s)$, for $s>p^n$, is an integral curve of degree $s^2+p^{2n}$ not lying on any surface of degree $s$ such that its generic plane section is contained in an integral plane curve of degree $s$.

\

I wish to express my deepest gratitude and appreciation to my Ph.D. advisor, Professor Rosario Strano, for his support and encouragement. My profound thanks and gratitude are addressed also to Riccardo Re for the many stimulating conversations and discussions.

\section{The Frobenius morphism}

First let us recall the definition of the relative Frobenius morphism (we follow Ein's notation in \cite{Ein}):
\begin{defin}
The absolute Frobenius morphism of a scheme $X$ of characteristic $p>0$ is $F_X\colon X\rightarrow X$, where $F_X$ is the identity as a map of topological spaces and on each $U$ open set $F^{\#}_{X}\colon\mathscr O_X(U)\rightarrow \mathscr O_X(U)$ is given by $f\mapsto f^p$ for each $f\in \mathscr O_X(U)$. Given $X\rightarrow S$ for some scheme $S$ and $X^{p/S}=X\times_{S,\, F_S}S$, the absolute Frobenius morphisms on $X$ and $S$ induce a morphism $F_{X/S}\colon X\rightarrow X^{p/S}$, called the Frobenius morphism of $X$ relative to $S$.
\end{defin}

Let now $r\in \mathbb N$ and $n\in \mathbb Z$ be integers. Let $F\colon \mathbb P^3\rightarrow \mathbb P^3$ be the absolute Frobenius and let us consider the sheaf $\mathscr F=(F^r)^{\star}\left(\Omega_{\mathbb P^3}\right)$. The following result will be needed later:
\begin{lemma}  \label{L:8}
\
\begin{enumerate}[label=\upshape(\roman*), labelindent=\parindent]
\item $h^0\left(\mathscr F(2p^r)\right)=6$,  \label{0}
\item $h^0\left(\mathscr F(n)\right)\ne 0$ if and only if $n\ge 2p^r$,   \label{1}
\item $h^2\left(\mathscr F(n)\right)=0$ for every $n\in \mathbb Z$, \label{2}
\end{enumerate}
\end{lemma}
\begin{proof}
First let us make some remarks. The sheaf $\Omega_{\mathbb P^3}$ is determined by the Euler sequence $0\rightarrow \Omega_{\mathbb P^3}\rightarrow \mathscr O^{\oplus 4}_{\mathbb P^3}(-1)\rightarrow \mathscr O_{\mathbb P^3}\rightarrow 0$, which is part of the Koszul complex $0 \rightarrow \mathscr O_{\mathbb P^3}(-4)\rightarrow \mathscr O^{\oplus 4}_{\mathbb P^3}(-3)\rightarrow \mathscr O^{\oplus 6}_{\mathbb P^3}(-2)\rightarrow \mathscr O^{\oplus 4}_{\mathbb P^3}(-1)\rightarrow\mathscr O_{\mathbb P^3}\rightarrow 0$. So $\mathscr F$, by the flatness of the absolute Frobenius, is determined by the exact sequence:
\begin{equation} \label{S:24}
0\rightarrow \mathscr F\rightarrow \mathscr O^{\oplus 4}_{\mathbb P^3}(-p^r)\rightarrow \mathscr O_{\mathbb P^3}\rightarrow 0,
\end{equation}
which is part of the following long exact sequence:
\begin{equation}  \label{S:25}
0 \rightarrow \mathscr O_{\mathbb P^3}(-4p^r)\rightarrow \mathscr O^{\oplus 4}_{\mathbb P^3}(-3p^r)\rightarrow \mathscr O^{\oplus 6}_{\mathbb P^3}(-2p^r)\rightarrow \mathscr O^{\oplus 4}_{\mathbb P^3}(-p^r)\rightarrow\mathscr O_{\mathbb P^3}\rightarrow 0.
\end{equation}
\ref{2} follows immediately from \eqref{S:24}. Now we prove \ref{0} and \ref{1}. Considered the cokernel $\mathscr G$ of the first nonzero map in \eqref{S:25}:
\begin{equation}  \label{S:26}
0 \rightarrow \mathscr O_{\mathbb P^3}(-4p^r)\rightarrow \mathscr O^{\oplus 4}_{\mathbb P^3}(-3p^r)\rightarrow \mathscr G\rightarrow 0
\end{equation}
$\mathscr F$ and $\mathscr G$ are related by the exact sequence:
\begin{equation}  \label{S:27}
0\rightarrow\mathscr G\rightarrow \mathscr O^{\oplus 6}_{\mathbb P^3}(-2p^r)\rightarrow \mathscr F\rightarrow 0.
\end{equation}
Since $\mathscr F$ and $\mathscr G$ are vector bundles and ${\mathscr F}^{\vee}\cong \mathscr G(4p^r)$, then they are reflexive and ${\mathscr G}^{\vee}\cong\mathscr F(4p^r)$. So by \eqref{S:26} and \eqref{S:27}:
\[h^0\left(\mathscr F(n)\right)=6h^0\left(\mathscr O_{\mathbb P^3}(n-2p^r)\right)-4h^0\left(\mathscr O_{\mathbb P^3}(n-3p^r)\right)+h^0\left(\mathscr O_{\mathbb P^3}(n-4p^r)\right).\]
From this we get \ref{0} and \ref{1}.
\end{proof}

In the notation of Lemma \ref{L:8}, let us consider the sheaf $\mathscr K=\mathscr F(p^r)|_H$,
restriction of $\mathscr F(p^r)$ to a plane $H$ in $\mathbb P^3$.
\begin{lemma}   \label{L:9}
For every $m\in \mathbb Z$:
\[h^0\left(\mathscr K(m)\right)=h^0\left(\mathscr O_H(m)\right)+3h^0\left(\mathscr O_H(m-p^r)\right)-h^0\left(\mathscr O_H(m-2p^r)\right).\]
\end{lemma}
\begin{proof}
Let us make the position $\mathscr F_H=(F^r)^{\star}(\Omega_H)$. We can construct a surjective morphism of sheaves $\varphi:{\mathscr O}^{\oplus 4}_H\rightarrow {\mathscr O}^{\oplus 3}_H$ in such a way that we get the following commutative diagram:
\[
\xymatrix@+1pc{
0 \ar[r] & \mathscr K \ar[d] \ar[r] & {\mathscr O}^{\oplus 4}_H \ar[d]_{\varphi} \ar[r] & \mathscr O_H(p^r) \ar[d]_{\text{id}} \ar[r] & 0\\
0 \ar[r] & \mathscr F_H(p^r) \ar[r] & {\mathscr O}^{\oplus 3}_H \ar[r] & \mathscr O_H(p^r) \ar[r] & 0.}
\]
Indeed, if $k[x_0,x_1,x_2,x_3]$ is the coordinate ring associated to $\mathbb P^3$ and $H$ has equation $x_3=\sum_{i=0}^2a_ix_i$, then we can define $\varphi$ as given by $\mathscr O^{\oplus 4}_{H}\ni(\sigma_0,\sigma_1,\sigma_2,\sigma_3)\mapsto (\sigma_0+{a_0}^{p^r}\sigma_3,\sigma_1+{a_1}^{p^r}\sigma_3,\sigma_2+{a_2}^{p^r}\sigma_3)\in \mathscr O^{\oplus 3}_{H}$.

So by the snake lemma and by the fact that $\operatorname{Ker}\varphi\cong \mathscr O_H$ we find the exact sequence:
\begin{equation} \label{S:28}
0\rightarrow \mathscr O_H\rightarrow \mathscr K\rightarrow \mathscr F_H(p^r)\rightarrow 0.
\end{equation}
Proceeding as in Lemma \ref{L:8} we see that $\mathscr F_H(p^r)$ comes from the Koszul complex $0\rightarrow \mathscr O_H(-2p^r)\rightarrow {\mathscr O}^{\oplus 3}_H(-p^r)\rightarrow {\mathscr O}^{\oplus 3}_H\rightarrow \mathscr O_H(p^r)\rightarrow 0$, which implies that $\mathscr F_H={\mathscr F_H}^{\vee}(-3p^r)$. Now:
\begin{multline} \label{S:29}
h^0\left(\mathscr F_H(p^r+m)\right)=h^0\left({\mathscr F_H}^{\vee}(m-2p^r)\right)=h^0\left((\mathscr F_H(2p^r-m))^{\vee}\right)=\\
=h^2\left(\mathscr F_H(2p^r-m)\otimes \mathscr O_H(-3)\right)=h^2\left(\mathscr F_H(2p^r-m-3)\right).
\end{multline}
From $0\rightarrow \mathscr F_H(2p^r-m-3)\rightarrow {\mathscr O}^{\oplus 3}_H(p^r-m-3)\rightarrow \mathscr O_H(2p^r-m-3)\rightarrow 0$
we see that $h^2\left(\mathscr F_H(2p^r-m-3)\right)=-h^0\left(\mathscr O_H(m-2p^r)\right)+3h^0\left(\mathscr O_H(m-p^r)\right)$, that, together with \eqref{S:28} and \eqref{S:29}, leads us to the conclusion.
\end{proof}

\section{Incidence varieties in characteristic $p$}

Let us consider the bi-projective space $\check{\mathbb P}^3\times\mathbb P^3$ and let $r\in \mathbb N$ be a non negative integer. Let $k[\underline t]$ and $k[\underline x]$ be the coordinate rings for $\check{\mathbb P}^3$ and $\mathbb P^3$, respectively. Let $M_r$ be the hypersurface of equation:
\[h_r:=\sum_{i=0}^3t_i{x_i}^{p^r}=0.\]
First we need the following result:
\begin{lemma}  \label{L:10}
Let $q\in k[\underline t,\underline x]$ be a homogeneous polynomial of bi-degree  $(\alpha,s)$ such that:
\begin{equation}  \label{S:30}
{x_i}^{p^r}\frac{\partial q}{\partial t_j}-{x_j}^{p^r}\frac{\partial q}{\partial t_i}\in (h_r)\ \forall\, i,j
\end{equation}
and $q\notin (h_r)$. Then there exists $q'=q+h_rm$ bi-homogeneous of bi-degree $(\alpha,s)$ such that:
\[\frac{\partial q'}{\partial t_i}=0\ \forall\, i.\]
\end{lemma}

Since the proof of this lemma requires some computations, we leave it to the end of this section. Let us now remark that in the case $r=0$ $M_r$ is usual incidence variety $M$ of equation $\sum t_ix_i=0$. If $r\ge 1$, $M_r$ is determined by the following fibred product:
\begin{equation} \label{S:35}
\xymatrix@+1pc{
M \ar@/_/[ddr]_{p} \ar@/^/[drr]^{(F_M)^r}
\ar@{-->}[dr]|-{F_{M_r}} \\
& M_r \ar[d]_(0.3){p_{M_r}} \ar[r]^{\pi} & M \ar[d]^{p} \\
& \mathbb P^3 \ar[r]_{F^r} & \mathbb P^3}
\end{equation}
where $F\colon \mathbb P^3\rightarrow\mathbb P^3$ is the absolute Frobenius. Moreover $M=\mathbb P(\Theta_{\mathbb P^3}(-1))$ and so by \cite[Lemma 1.5]{Ein} we get $M_r=\mathbb P({F^r}^{\star}(\Theta_{\mathbb P^3}(-1)))$. By \cite[Ch.II, ex. 7.8]{Hart} this implies:
\begin{equation} \label{pic}
\operatorname{Pic}(M_r)=\mathbb Z\times \mathbb Z
\end{equation}
for any $r\ge 0$.

Let us consider an integral hypersurface $V\subset M_r$ and let us suppose that the projection $\pi\colon V\rightarrow \mathbb P^3$ is dominant. Using the previous lemma we prove the following result:
\begin{prop}  \label{P:14}
If $\pi$ is not generically smooth, then there exists $s\ge 1$, such that $V\subset \check{\mathbb P}^3\times \mathbb P^3$ is the complete intersection determined by $g=h_r=0$, for some $g\in k[{\underline t}^{p^s},\underline x]$.
\end{prop}
\begin{proof}
Since $M_r\subset \check{\mathbb P}^3\times \mathbb P^3$ is a hypersurface of bi-degree $(1,p^r)$, the structure sheaf $\mathscr O_{M_r}$ is given by $0\rightarrow \mathscr O_{\check{\mathbb P}^3\times \mathbb P^3}(-1,-p^r)\rightarrow \mathscr O_{\check{\mathbb P}^3\times \mathbb P^3}\rightarrow \mathscr O_{M_r}\rightarrow 0$. By the K\"unneth formula (\cite[Ch.~VI, Corollary~8.13]{Milne80}) $h^1\left(\mathscr O_{\check{\mathbb P}^3\times \mathbb P^3}(m,n)\right)=0$ for every $m$, $n\in \mathbb Z$, so that the morphism $H^0\left(\mathscr O_{\check{\mathbb P}^3\times \mathbb P^3}(m,n)\right)\rightarrow H^0\left(\mathscr O_{M_r}(m,n)\right)$ is surjective for every $m$, $n\in \mathbb Z$. Together with \eqref{pic} this implies that $V\subset \check{\mathbb P}^3\times \mathbb P^3$ is a complete intersection given by $g=h_r=0$ for some $g\in k[\underline t,\underline x]$ bi-homogeneous of bi-degree $(m,n)$ for some $m$, $n\in \mathbb N$.

Let $P_0=(\underline a,\underline b)\in V$ be a regular point. By hypothesis the map on the projective tangent spaces $T_{V,P_0}$ and $T_{\mathbb P^3\!,{\pi}(P_0)}$ is not surjective. The projective tangent space $T_{V,P_0}$ at $P_0\in V$ is given by the equations:
\[
\sum_{i=0}^3 \frac{\partial g}{\partial x_i}(P_0)x_i+\sum_{i=0}^3\frac{\partial g}{\partial t_i}(P_0)t_i=\sum_{i=0}^3(a_ix_i+b_it_i)=0\]
if $r=0$ and by the equations:
\[
\sum_{i=0}^3 \frac{\partial g}{\partial x_i}(P_0)x_i+\sum_{i=0}^3\frac{\partial g}{\partial t_i}(P_0)t_i=\sum_{i=0}^3{b_i}^{p^r}t_i=0
\]
if $r\ge 1$. In both cases the projection on $T_{\mathbb P^3\!,{\pi}(P_0)}$ is not surjective if and only if there exists $\lambda\in k$ such that:
\[\dfrac{\partial g}{\partial t_i}(P_0)=\lambda {b_i}^{p^r}\quad \forall\, i=0,\dots,3.\]
So in such a situation:
\[{b_i}^{p^r}\frac{\partial g}{\partial t_j}(P_0)-{b_j}^{p^r}\frac{\partial g}{\partial t_i}(P_0)=0\quad \forall\, i,j.\]
This means that for every $i$, $j$ the hypersurface $V_{ij}\subset \check{\mathbb P}^3\times \mathbb P^3$
given by ${x_i}^{p^r}\tfrac{\partial g}{\partial t_j}-{x_j}^{p^r}\tfrac{\partial g}{\partial t_i}=0$ contains Reg\,$(V)$, the open subset of the regular points of $V$. So $V_{ij}\supset V$ for all $i$, $j$, which means that:
\[{x_i}^{p^r}\frac{\partial g}{\partial t_j}-{x_j}^{p^r}\frac{\partial g}{\partial t_i}\in (g,h_r)\quad \forall\, i,j.\]
If ${x_i}^{p^r}\tfrac{\partial g}{\partial t_j}-{x_j}^{p^r}\tfrac{\partial g}{\partial t_i}$ is a nonzero polynomial, then it is a bi-homogeneous polynomial of bi-degree $(m-1,n+p^r)$. Since $g$ is bi-homogeneous of bi-degree $(m,n)$, then
\[{x_i}^{p^r}\frac{\partial g}{\partial t_j}-{x_j}^{p^r}\frac{\partial g}{\partial t_i}\in (h_r)\quad \forall\, i,j.\]
Applying Lemma \ref{L:10} we see that there exists $m\in k[\underline t,\underline x]$ such that, given $g'=g+mh_r$, we have ${\partial g'}/{\partial t_i}=0$ for every $i$. So by replacing $g$ by $g'$ we can suppose that $g\in k[{\underline t}^{p^s},\underline x]$, for some $s\ge 1$.
\end{proof}

Now we can prove the main result of this section:
\begin{thm}  \label{T:21}
Let $V\subset \check{\mathbb P}^3\times \mathbb P^3$ be an integral hypersurface in $M$ such that the projection $\pi\colon V\rightarrow \mathbb P^3$ is dominant and not generically smooth. Then there exist $r\ge 1$, and $V_r\subset M_r$ integral hypersurface such that $\pi$ can be factored in the following way:
\[
\xymatrix@+1pc{
V \ar[rr]^{\pi} \ar[dr]_{F_r}& & \mathbb P^3\\
& V_r \ar[ur]_{\pi_r}
}
\]
where the projection $\pi_r$ is dominant and generically smooth and $F_r$ is induced by the commutative diagram:
\[
\xymatrix@+1pc{
V \ar[r]^{F_r} \ar[d]_{j} & V_r \ar[d]^{i}\\
M \ar[r]_{F_{M_r}} & M_r.}
\]
\end{thm}
\begin{proof}
First note that by hypothesis and by Proposition \ref{P:14} it follows that $V\subset \check{\mathbb P}^3\times \mathbb P^3$ is the complete intersection determined by $h=0$ and $q=0$, for some $q\in k[{\underline t}^{p^r},\underline x]$ and $r\ge 1$. We can suppose that $q\in k[{\underline t}^{p^r},\underline x]$ and $q'\notin k[{\underline t}^{p^{r+1}},\underline x]$ for any $q'\equiv q$ $\bmod\, (h)$. So we can say that $q(\underline t,\underline x)=f({\underline t}^{p^r},\underline x)$ for some bi-homogeneous $f\in k[\underline t,\underline x]$.

Let us now consider the hypersurface $M_r\subset \check{\mathbb P}^3\times \mathbb P^3$ and the two projections $p_{M_r}\colon M_r\rightarrow \mathbb P^3$ and $g_{M_r}\colon M_r\rightarrow \check{\mathbb P}^3$. Let $V_r\subset M_r$ be the hypersurface determined by $f=0$. Since the morphism $F_{M_r}\colon M\rightarrow M_r$ in diagram \eqref{S:35} is the $p^r$-th power on the $\{\mspace{1mu}t_i\mspace{1mu}\}$, we see that $V$ is the following fibred product:
\begin{equation}  \label{S:36}
\xymatrix@+1pc{
V \ar[r]^{F_r} \ar[d]_{j} & V_r \ar[d]^{i}\\
M \ar[r]_{F_{M_r}} & M_r.}
\end{equation}
So, being $V$ integral, $V_r$ must be integral too.

$\pi_r$ is dominant because the composition $V\rightarrow V_r\rightarrow \mathbb P^3$ is dominant. Now we show that $\pi_r$ is generically smooth. $V_r\subset\check{\mathbb P}^3\times \mathbb P^3$ is determined by the complete intersection $\sum_{i=0}^3 t_i{x_i}^{p^r}=f(\underline t,\underline x)=0$. If $\pi_r$ were not generically smooth, then by Proposition \ref{P:14} we could suppose that $f\in k[{\underline t}^p,\underline x]$. This would imply, by the commutative diagram \eqref{S:36}, that we could take $q\in k[{\underline t}^{p^{r+1}},\underline x]$, which contradicts the choice of $r$.
\end{proof}

Let us now return to Lemma \ref{L:10}.
\begin{proof}[Proof of Lemma \ref{L:10}]
From \eqref{S:30} we have:
\[\left(\sum_{i=0}^3t_i{x_i}^{p^r}\right)\frac{\partial q}{\partial t_j}-{x_j}^{p^r}\left(\sum_{i=0}^3t_i\frac{\partial q}{\partial t_i}\right)\in (h_r)\ \forall\, j.\]
Using that $\sum t_i{x_i}^{p^r}=h_r$ and that $h_r$ is irreducible, we deduce $\sum t_i\partial q/\partial t_i\in (h_r)$. However $\sum_{i=0}^3t_i\frac{\partial q}{\partial t_i}=aq$, where $a$ is the remainder of the division of $\alpha$ by $p$, because $q$ is homogeneous of degree $\alpha$ in the $\{\mspace{1mu}t_i\mspace{1mu}\}$. So $aq\in (h_r)$ and by hypothesis the only possibility is that $a=0$, which means that $p\mid\alpha$.

By \eqref{S:30} for every $i$, $j=0$, $1$, $2$, $3$ there exists $l_{ij}$ bi-homogeneous in $k[\underline t,\underline x]$ such that:
\begin{equation} \label{S:31}
{x_i}^{p^r}\frac{\partial q}{\partial t_j}-{x_j}^{p^r}\frac{\partial q}{\partial t_i}=l_{ij}h_r.
\end{equation}
The identity:
\begin{multline*}
{x_k}^{p^r}\left({x_i}^{p^r}\frac{\partial q}{\partial t_j}-{x_j}^{p^r}\frac{\partial q}{\partial t_i}\right)-{x_i}^{p^r}\left({x_k}^{p^r}\frac{\partial q}{\partial t_j}-{x_j}^{p^r}\frac{\partial q}{\partial t_k}\right)+{}\\
{}+{x_j}^{p^r}\left({x_k}^{p^r}\frac{\partial q}{\partial t_i}-{x_i}^{p^r}\frac{\partial q}{\partial t_k}\right)=0
\end{multline*}
for every $i,j,k$ determines the equality ${x_k}^{p^r}l_{ij}-{x_i}^{p^r}l_{kj}+{x_j}^{p^r}l_{ki}=0$. So on $D_+(x_ix_jx_k)$ we have the equality:
\[\frac{l_{ij}}{(x_ix_j)^{p^r}}-\frac{l_{ki}}{(x_kx_j)^{p^r}}+\frac{l_{ki}}{(x_kx_i)^{p^r}}=0.\]
Considered now the open covering $\mathfrak U=\{D_+(x_i)\mid i=0,\dots,3\}$ and $n=\deg l_{ij}$, we get a \v{C}ech cocycle in $\check{H}^1(\mathfrak U, \mathscr O_{\mathbb P^3}(n-2p^r))\cong H^1(\mathscr O_{\mathbb P^3}(n-2p^r))=0$. So the cocycle is a coboundary and for every $i$, $j$ there exist $m_i$, $m_j\in k[\underline t,\underline x]$ such that $l_{ij}=m_i{x_j}^{p^r}-m_j{x_i}^{p^r}$. So by \eqref{S:31}:
\[{x_i}^{p^r}\left(\frac{\partial q}{\partial t_j}-m_jh_r\right)={x_j}^{p^r}\left(\frac{\partial q}{\partial t_i}-m_ih_r\right)\ \forall\, i,j.\]
So there exists $m\in k[\underline t,\underline x]$ such that $\partial q/\partial t_i=m_ih_r+m{x_i}^{p^r}$ for every $i$. By replacing $q$ by $q-mh_r$ we may assume that ${\partial q}/{\partial t_i}=m_ih_r$ for every $i$. We want to prove that $\partial q /\partial t_i=U_i{h_r}^p$ for some $U_i$ and so let us suppose that:
\[\frac{\partial q}{\partial t_i}=v_i{h_r}^n\ \forall\, i\]
for some $n<p-1$. Then:
\[\frac{{\partial}^2 q}{\partial t_it_j}=\frac{\partial v_i}{\partial t_j}{h_r}^n+nv_i{x_j}^{p^r}{h_r}^{n-1}\ \forall\, i,j.\]
But we have also:
\[\frac{{\partial}^2 q}{\partial t_it_j}=\frac{\partial v_j}{\partial t_i}{h_r}^n+nv_j{x_i}^{p^r}{h_r}^{n-1}\ \forall\, i,j.\]
So:
\[\frac{\partial v_i}{\partial t_j}{h_r}^n+nv_i{x_j}^{p^r}{h_r}^{n-1}=\frac{\partial v_j}{\partial t_i}{h_r}^n+nv_j{x_i}^{p^r}{h_r}^{n-1}\]
\[\Rightarrow v_i{x_j}^{p^r}-v_j{x_i}^{p^r}=\frac{h_r}{n}\left(\frac{\partial v_j}{\partial t_i}-\frac{\partial v_i}{\partial t_j}\right)\ \forall\, i,j.\]
This implies that $v_i=v{x_i}^{p^r}+h_ru_i$ for every $i$. By replacing $q$ by $q-\tfrac{1}{n+1}v{h_r}^{n+1}$ we may assume that:
\[\frac{\partial q}{\partial t_i}=V_i{h_r}^{p-1}\ \forall\, i.\]
We know that:
\begin{equation} \label{S:33}
\frac{{\partial}^{p}q}{\partial {t_i}^p}=0.
\end{equation}
This means that:
\[\frac{{\partial}^{p-1}(V_i{h_r}^{p-1})}{\partial {t_i}^{p-1}}=0\]
\[\Rightarrow \sum_{n=0}^{p-1}\binom{p-1}{n} \frac{{\partial}^n V_i}{\partial {t_i}^{n}}\frac{{\partial}^{p-1-n}({h_r}^{p-1})}{\partial {t_i}^{p-1-n}}=0\]
\[\Rightarrow h_r\mid (p-1)!{x_i}^{p^{r+1}-p^r}V_i\quad\Rightarrow h_r\mid V_i\ \forall\, i.\]
So we can suppose that:
\[\frac{\partial q}{\partial t_i}=U_i{h_r}^p\ \forall\, i.\]
Now \eqref{S:33} leads us to the conclusion that:
\[\frac{{\partial}^{p-1} U_i}{\partial {t_i}^{p-1}}=0\]
which means that in $U_i$, for each $i$, there are no terms of type ${t_i}^{kp-1}$ for any $k\ge 1$. So in particular we can say that:
\[U_0=\frac{\partial M_0}{\partial t_0}\]
for some $M_0$ bi-homogeneous. Now $q'=q-M_0{h_r}^p$ is such that:
\[\frac{\partial q'}{\partial t_0}=0\quad\mbox{and}\quad\frac{\partial q'}{\partial t_i}=U'_i{h_r}^p,\ i=1,2,3\]
\[\Rightarrow \frac{\partial U'_i}{\partial t_0}=0,\ i=1,2,3.\]
So we can find $U''_1$ such that:
\[\frac{\partial U''_1}{\partial t_0}=0\quad\mbox{and}\quad\frac{\partial U''_1}{\partial t_1}=U'_1.\]
If we consider $q''=q'-U''_1{h_r}^p$ we see that:
\[\frac{\partial q''}{\partial t_i}=0,\ i=0,1\quad\mbox{and}\quad\frac{\partial q''}{\partial t_i}=U''_i{h_r}^p,\ i=2,3.\]
Proceeding in this way we get $\partial q/{\partial t_i}=0$ for every $i$.
\end{proof}

\section{Proof of the main theorem}

Let us consider now a curve $C\subset \mathbb P^3$ and, following the notation of Theorem \ref{T:21}, the projections $p_{M_r}\colon M_r\rightarrow \mathbb P^3$ and $g_{M_r}\colon M_r\rightarrow \check{\mathbb P}^3$. Let $T_r={p_{M_r}}^{-1}(C)$ and:
\[\mathscr{I}_r(m,n)={g_{M_r}}^{\star}\left(\mathscr O_{{\mathbb P^3}^{\vee}}(m)\right)\otimes_{\mathscr O_{M_r}}{p_{M_r}}^{\star}\left(\mathscr I_C(n)\right)\]
for every $m$, $n\in \mathbb Z$.
\begin{prop}  \label{P:15}
If $\mathscr I_r=\mathscr I_r(0,0)$ and $\mathscr I_{T_r}$ is the ideal sheaf of $T_r$ in $M_r$, then $\mathscr I_r=\mathscr I_{T_r}$.
\end{prop}
\begin{proof}
First note that $\mathscr I_r(m,n)=\mathscr O_{M_r}(m,n)\otimes_{\mathscr O_{M_r}}{p_{M_r}}^{\star}\left(\mathscr I_C\right)$ for any $m,n\in \mathbb Z$. Moreover $p_M$ is smooth, in particular flat. So, by base change (see \eqref{S:35}), $p_{M_r}$ is flat too and we can apply \cite[Ch.~III, Proposition 9.3]{Hart} to the following commutative diagram:
\[
\xymatrix@+1pc{
T_r \ar[r]^{\pi} \ar[d]_{j} & C \ar[d]^{i}\\
M_r \ar[r]_{p_{M_r}} & \mathbb P^3}
\]
to get that ${p_{M_r}}^{\star}i_{\star}\mathscr O_C\cong j_{\star}{\pi}^{\star}\mathscr O_C\cong j_{\star}\mathscr O_{T_r}$. This fact together with the exact sequence $0 \rightarrow {p_{M_r}}^{\star} \mathscr I_C \rightarrow {p_{M_r}}^{\star}\mathscr O_{\mathbb P^3} \rightarrow {p_{M_r}}^{\star}i_{\star}\mathscr O_C \rightarrow 0$, consequence of the flatness of $p_{M_r}$, leads us to the desired conclusion.
\end{proof}

Now we are ready to prove Theorem \ref{T:1}.

\begin{proof}[Proof of Theorem \ref{T:1}]
We divide the proof in different steps.
\begin{step}
There exists $S\subset \check{\mathbb P}^3\times \mathbb P^3$ integral such that the generic fibre of the projection $S\rightarrow \check{\mathbb P}^3$ is an integral plane curve of degree $s$ containing $\Gamma$.
\end{step}
\begin{proof}[Proof of Step 1]
Let $\mathscr I_C$ be the ideal sheaf of $C$ in $\mathbb P^3$ and let $M\subset \check{\mathbb P}^3\times\mathbb P^3$ be the incidence variety. Let us consider the two projections:
\[
\xymatrix@+0.5pc{
& M \ar[dl]_{p} \ar[dr]^{g} &\\
\mathbb P^3 & & \check{\mathbb P}^3
}
\]
and the $\mathscr O_M$-module $\mathscr I(m,n)={g}^{\star}\left(\mathscr O_{{\mathbb P^3}^{\vee}}(m)\right)\otimes_{\mathscr O_M}{p}^{\star}\left(\mathscr I_C(n)\right)$.
As we have seen in Proposition \ref{P:15} in the case $r=0$ $\mathscr I$ is the ideal sheaf of $T={p}^{-1}(C)$ in $\mathscr O_M$. Moreover there exists $\alpha$ such that $H^0(\mathscr I(\alpha,s))\ne 0$. Indeed, if $\eta\in\check{\mathbb P}^3$ denotes the generic point and $\Gamma$ is the generic plane section of $C$, then $H^0(p^{\star}\mathscr I(s)|_{g^{-1}(\eta)})=H^0(\mathscr I_{\Gamma}(s))\ne 0$ and so this global section determines an effective divisor in $M_{k(\eta)}=M\times_{\check{\mathbb P}^3} \operatorname{Spec}k(\eta)\cong \mathbb P^2_{k(\eta)}$. Then there exists $U\subset \check{\mathbb P}^3$ such that this divisor extends to an effective divisor $D\subset M_U=M\times_{\check{\mathbb P}^3}U$ containing $T\times_{\check{\mathbb P}^3} U$. The closure $\overline D\subset M$ is an effective divisor containing $T$ and, since $\operatorname{Pic}(M)=\mathbb Z\times \mathbb Z$, it is a divisor determined by a global section of $\mathscr I(\alpha,s)$, for some $\alpha\ge 0$.

Taken the least $\alpha$ such that $h^0(\mathscr I(\alpha,s))\ne 0$, there exists $q\in H^0(\mathscr I(\alpha,s))$ that determines a hypersurface $S$ in $M$ such that $S\cap {g}^{-1}(\eta)$ is an integral curve of degree $s$ containing $\Gamma$. Moreover, as we saw in Proposition \ref{P:14}, $S$ is a complete intersection of codimension $2$ in $\check{\mathbb P}^3\times\mathbb P^3$. This implies that $S$ is irreducible.
\renewcommand{\qedsymbol}{$\Box$}
\end{proof}

To prove Theorem \ref{T:1}, we now assume that the curve $C$ is not contained in any surface of degree $s$, in other words, $h^0(\mathscr I_C(s))=0$. Then $p_S\colon S\rightarrow \mathbb P^3$ is dominant and, since $\alpha\ge 0$, in such a situation it must be $\alpha>0$.

\begin{step}
We can factor $p_S$ through a generically smooth morphism $S_r\rightarrow \mathbb P^3$, with $S_r$ scheme of zeroes of a global section of $\mathscr I_r(\beta,s)$, being $\mathscr I_r={p_{M_r}}^{\star} \mathscr I_C$, and $\alpha=\beta p^r$, for some $r\ge 0$.
\end{step}
\begin{proof}[Proof of Step 2]
If $p_S$ is not generically smooth, then by Theorem \ref{T:21} it follows that there exist $r\ge 1$, and $S_r\subset M_r$ integral hypersurface such that $p_S$ can be factored in the following way:
\[
\xymatrix@+1pc{
S \ar[rr]^{p_S} \ar[dr]_{F^r}& & \mathbb P^3\\
& S_r \ar[ur]_{p_{S_r}}
}
\]
where the projection $p_{S_r}$ is dominant and generically smooth and $F^r$ is induced by the commutative diagram:
\[
\xymatrix@+1pc{
S \ar[r]^{F^r} \ar[d]_{j} & S_r \ar[d]^{i}\\
M \ar[r]_{F_{M_r}} & M_r.}
\]
Moreover, we also get that $\alpha=p^r\beta$, for some $\beta\in \mathbb N$, $\beta>0$.

Considered the sheaf $\mathscr{I}_r={p_{M_r}}^{\star}\mathscr I_C$ and the scheme $T_r={p_{M_r}}^{-1}(C)$, by Proposition \ref{P:15} we see that $\mathscr I_r$ is the ideal sheaf of $T_r$ in $M_r$. Given $T={p}^{-1}(C)$, since $S\supset T$ and $F_{M_r}(T)=T_r$, then $S_r\supset T_r$. So $S_r\subset M_r$ is the scheme of zeros of a global section in $H^0(\mathscr{I}_r(\beta,s))$.

Hence in both cases we find $S_r$ integral hypersurface in $M_r$, with $r\ge 0$, such that the projection $p_{S_r}\colon S_r\rightarrow \mathbb P^3$ is generically smooth and $S_r\subset M_r$ is the scheme of zeros of a global section in $H^0(\mathscr{I}_r(\beta,s))$, for some $\beta>0$.
\renewcommand{\qedsymbol}{$\Box$}
\end{proof}

Let us now follow the proof of Gruson and Peskine given in \cite{GruPes}.

\begin{step}
There exists a 3-dimensional scheme $Y$, with $T_r\subseteq Y\subset S_r$, such that we have:
\begin{equation}   \label{S:1}
0\rightarrow \Omega^{\vee}_{S_r/\mathbb P^3}\rightarrow \Omega^{\vee}_{M_r/\mathbb P^3}\otimes_{\mathscr O_{M_r}}\mathscr O_{S_r}\rightarrow \mathscr I_{Y}(\beta,s)\rightarrow 0
\end{equation}
with $\mathscr I_Y\subset \mathscr O_{S_r}$ ideal sheaf of $Y$.
\end{step}
\begin{proof}[Proof of Step 3]
Since $S_r$ is generically smooth over $\mathbb P^3$, we have the exact sequence $0\rightarrow {\mathscr O_{S_r}(-\beta,-s)}\rightarrow \Omega_{M_r/{\mathbb P^3}}\otimes_{\mathscr O_{M_r}} \mathscr O_{S_r}\rightarrow \Omega_{S_r/\mathbb P^3}\rightarrow 0$. Dualizing with respect to $\mathscr O_{S_r}$, we get:
\begin{equation} \label{S:37}
0\rightarrow \Omega^{\vee}_{S_r/\mathbb P^3}\rightarrow \Omega^{\vee}_{M_r/\mathbb P^3}\otimes_{\mathscr O_{M_r}}\mathscr O_{S_r}\rightarrow \mathscr O_{S_r}(\beta,s).
\end{equation}
Since all the fibres of the projection $T_r\rightarrow C$ have dimension 2 and $\dim S_r=4$, $p_{S_r}$ is not regular in any of the points of $T_r$. It means that the last map in \eqref{S:37} has image inside $\mathscr I_{T_r}(\beta,s)$, the ideal sheaf of $T_r$ in $S_r$. So this image is an ideal sheaf of type $\mathscr I_Y(\beta,s)$, where $Y\subset S_r$ is a scheme containing $T_r$, and $3=\dim T_r\le \dim Y\le \dim S_r=4$. Since $S_r$ is reduced and irreducible, if $\dim Y=4$, then $p_{S_r}$ would be non regular almost everywhere in $S_r$. This contradicts the fact that $p_{S_r}$ is generically smooth. So $\dim Y=\dim T_r=3$ and $T_r\subseteq Y$.
\renewcommand{\qedsymbol}{$\Box$}
\end{proof}

Let us consider the projection $g_{M_r}\colon M_r\rightarrow \check{\mathbb P}^3$ and take any $(\underline b)=({\underline d}^{p^r})\in \check{\mathbb P}^3$. Then ${g_{M_r}}^{-1}(\underline b)=\{(\underline x,{\underline d}^{p^r})\mid (\sum d_ix_i)^{p^r}=0\}$. If $H={{g_{M_r}}^{-1}(\underline b)}_{\text{red}}$ and $D=p(g^{-1}(\underline b)_{\text{red}})$, then there exists $U\subset\check{\mathbb P}^3$ open such that, taken $(\underline b)\in U$, $D$ is an irreducible curve of degree $s$ containing the plane section of $C$ with $H$. Let $\Gamma$ denote such a section and let $\mathscr I_{\Gamma}\subset \mathscr O_D$ be the its ideal sheaf.

\begin{step}
If $\mathscr M=(\Omega_{M_r/\mathbb P^3})|_H$, there exist a rank two vector bundle $\mathscr N$ and a zero-dimensional scheme $\Delta$, with $\Gamma\subseteq\Delta\subset D$, such that we have:
\begin{equation}  \label{S:40}
0\rightarrow \mathscr N\rightarrow \mathscr M^{\vee}\rightarrow \mathscr I_{\Delta}(s)\rightarrow 0
\end{equation}
being $\mathscr I_{\Delta}\subset \mathscr O_D$ the ideal sheaf of $\Delta$.
\end{step}
\begin{proof}[Proof of Step 4]
Since $M=\mathbb P(\Theta_{\mathbb P^3}(-1))$, by \eqref{S:35} and by \cite[Lemma 1.5]{Ein} we see that $M_r=\mathbb P({F^r}^{\star}(\Theta_{\mathbb P^3}(-1)))$, where we denoted by $F$, as in \eqref{S:35}, the absolute Frobenius on $\mathbb P^3$. The sheaf $\mathscr E={F^r}^{\star}(\Theta_{\mathbb P^3}(-1))$ is determined by the exact sequence $0\rightarrow \mathscr O_{\mathbb P^3}(-p^r)\rightarrow {\mathscr O_{\mathbb P3}}^{\oplus 4}\rightarrow \mathscr E\rightarrow 0$ and by \cite[Ch. III, Ex. 8.4(b)]{Hart} we have also $0\rightarrow \Omega_{M_r/\mathbb P^3}\rightarrow \big({p_{M_r}}^{\star}\mathscr E\big)(-1)\rightarrow \mathscr O_{M_r}\rightarrow 0$. When we restrict to $H$, by the fact that the sequences locally split it follows that the following sequences are exact:
\begin{equation} \label{S:38}
0\rightarrow \mathscr O_H(-p^r)\rightarrow {\mathscr O_H}^{\oplus 4}\rightarrow {\mathscr E}_H\rightarrow 0
\end{equation}
and
\[0\rightarrow \mathscr M\rightarrow \big({p_{M_r}}^{\star}\mathscr E\big)|_H(-1)\rightarrow \mathscr O_H\rightarrow 0\]
where $\mathscr M=(\Omega_{M_r/\mathbb P^3})|_H$. Since $({p_{M_r}}^{\star}\mathscr E\big)|_H(-1)=({p_{M_r}}^{\star}\mathscr E(-1,0)\big)|_H={\mathscr E}_H$, we have:
\begin{equation}   \label{S:39}
0\rightarrow \mathscr M\rightarrow \mathscr E_H\rightarrow \mathscr O_H\rightarrow 0.
\end{equation}

Restricting \eqref{S:1} to $H$, we get a surjective map ${\mathscr M}^{\vee}\otimes_{\mathscr O_H}\mathscr O_D\rightarrow \mathscr I_{\Delta}(s)$, with $\mathscr I_{\Delta}\subset \mathscr O_D$ ideal sheaf of a zero-dimensional scheme $\Delta$ containing $\Gamma$. The kernel of this map is a locally free sheaf of rank 2 that determines the exact sequence \eqref{S:40}.
\renewcommand{\qedsymbol}{$\Box$}
\end{proof}

\begin{step}
$d\le s^2+p^{2r}$.
\end{step}
\begin{proof}[Proof of Step 5]
Note that $c_1\big(\mathscr I_{\Delta}(s)\big)=s$ and $c_2\big(\mathscr I_{\Delta}(s)\big)=\deg \Delta=\delta\ge d$. Now we compute the Chern classes of the other sheaves. From \eqref{S:38} we have $c_1\big(\mathscr E_H\big)=p^r$ and $c_2\big(\mathscr E_H\big)=p^{2r}$. So by \eqref{S:39} $c_1(\mathscr M)=p^r$ and $c_2(\mathscr M)=p^{2r}$, from which it follows that $c_1\big({\mathscr M}^{\vee}\big)=-p^r$ and $c_2\big({\mathscr M}^{\vee}\big)=p^{2r}$. By \eqref{S:40} we see that:
\[c_1(\mathscr N)=-p^r-s\quad \mbox{and}\quad c_2(\mathscr N)=p^{2r}-\delta+s^2+p^rs.\]

Let $m\in \mathbb Z$ be the smallest number such that $H^0(\mathscr {\mathscr M}^{\vee}(m-1))=0$ and $H^0(\mathscr {\mathscr M}^{\vee}(m))>0$. Dualizing \eqref{S:39}, since $\mathscr O_H$ is a locally free sheaf, we get $0\rightarrow \mathscr O_H\rightarrow \mathscr E^{\vee}_H\rightarrow {\mathscr M}^{\vee}\rightarrow 0$, from which it follows that:
\[h^0\big({\mathscr M}^{\vee}(m)\big)=h^0\big(\mathscr E^{\vee}_H(m)\big)-h^0\big(\mathscr O_H(m)\big)\ \forall m\in \mathbb Z.\]
From the exact sequence \eqref{S:38} we see that, in the notation of Lemma \ref{L:9}, $\mathscr E_H^{\vee}=\mathscr K$, so that, by Lemma \ref{L:9} we see that:
\[h^0\left(\mathscr E_H^{\vee}(m)\right)=h^0\left(\mathscr O_H(m)\right)+3h^0\left(\mathscr O_H(m-p^r)\right)-h^0\left(\mathscr O_H(m-2p^r)\right)\]
for every $m\in \mathbb Z$. So:
\[h^0\big({\mathscr M}^{\vee}(m)\big)=3h^0\big(\mathscr O_H(m-p^r)\big)-h^0\big(\mathscr O_H(m-2p^r)\big)\ \forall m\in \mathbb Z.\]
This implies that $h^0({\mathscr M}^{\vee}(p^r-1))=0$ and $h^0({\mathscr M}^{\vee}(p^r))>0$. So $h^0(\mathscr N(p^r-1))=0$ and $p^{2r}+p^r(-s-p^r)+c_2\big(\mathscr N\big)=c_2(\mathscr N(p^r))\ge 0$ by \cite{GruPes}. So we get that $p^{2r}+s^2\ge \delta$ and, since $\delta \ge d$:
\begin{equation} \label{D:2}
p^{2r}+s^2\ge d.
\end{equation}
\renewcommand{\qedsymbol}{$\Box$}
\end{proof}

\begin{step}
If $C$ is connected, $p^r<s$; if $C$ is merely reduced, $p^r\le s$.
\end{step}
\begin{proof}[Proof of Step 6]
Let us now consider a generic plane $H=\operatorname{V}(l)$, with $l$ linear form in the $\{\mspace{1mu}x_i\mspace{1mu}\}$, and the non reduced surface $H_r$ in $\mathbb P^3$ given by $l^{p^r}=0$. Let $\Gamma_r\subset H_r$ be the section of $C$ with $H_r$. Then there is the following exact sequence:
\[0\rightarrow \mathscr I_C(-p^r)\stackrel{\varphi_H}{\rightarrow}\mathscr I_C\rightarrow i_{\star}\mathscr I_{\Gamma_r}\rightarrow 0\]
where $i\colon\Gamma_r\hookrightarrow \mathbb P^3$ and $\varphi_H$ is the multiplication by $l^{p^r}$. The long cohomology exact sequence associated to the previous exact sequence shifted by $s$ determines the following one:
\begin{equation} \label{S:41}
H^0\left(\mathscr I_C(s)\right)\rightarrow H^0\left(\mathscr I_{\Gamma_r}(s)\right)\rightarrow H^1\left(\mathscr I_C(s-p^r)\right)\stackrel{\varphi_H}{\rightarrow}H^1\left(\mathscr I_C(s)\right).
\end{equation}
Let $[H]\in \check{\mathbb P}^3$ be a point such that the fibre of the projection $M_r\rightarrow \check{\mathbb P}^3$ at $[H]$ is isomorphic to $H_r$. Then, taking $[H]$ in a suitable open $U\subset\check{\mathbb P}^3$, ${g_{S_r}}^{-1}([H])$ is the complete intersection of $H_r$ and of a surface of degree $s$ containing $C\cap H_r$, because $T_r\subset S_r$. It means that $H^0(\mathscr I_{\Gamma_r}(s))\ne 0$ and so by \eqref{S:41} and by hypothesis it must be $h^1(\mathscr I_C(s-p^r))\ne 0$.

Let us suppose that $C$ is connected. Then $h^1(\mathscr I_C(n))=0$ for $n\le 0$. So $s-p^r\ge 1$, because otherwise $h^0(\mathscr I_C(s))\ne 0$, which contradicts the hypothesis made at the beginning.

If $C$ is merely reduced, we still have $h^1(\mathscr I_C(n))=0$ for $n<0$. So, as before, it must be $s-p^r\ge 0$.
\renewcommand{\qedsymbol}{$\Box$}
\end{proof}

Let us suppose that $C$ is connected. By Step 6, if $p\ge s$, then the only possibility is that $r=0$, which implies $d\le s^2+1$. If $p<s$, then $p^r\le s-1$ and, in particular, $p^r\le p^n$, being $p^n<s\le p^{n+1}$. So by \eqref{D:2} $d\le s^2+p^{2n}$ and in particular we see that $d\le 2s^2-2s+1$.

Let us suppose now that $C$ is merely reduced. If $p>s$, then it must be $r=0$, so that by \eqref{D:2} we have $d\le s^2+1$. If $p\le s$, then $p^r\le s$ and, in particular, $p^r\le p^n$, being $p^n\le s<p^{n+1}$. Now by \eqref{D:2} we see that $d\le s^2+p^{2n}$. In particular, $d\le 2s^2$.
\end{proof}

\section{Example}

In this section we show that for any $p$ there exist smooth integral curves of degree $d=s^2+p^{2n}$, being $s>p$ and $n$ such that $p^n<s\le p^{n+1}$, that are not contained in any surface of degree $s$ and that have the generic plane section contained in an integral plane curve of degree $s$.

First, let us recall the following definition.

\begin{defin}
A rank $2$ vector bundle $\mathscr E_0$ on $\mathbb P^3$ is said to be a \emph{null-correlation bundle} if there exists an exact sequence:
\begin{equation}  \label{S:42}
0 \rightarrow \mathscr O_{\mathbb P^3} \stackrel{\tau}{\rightarrow} \Omega_{\mathbb P^3}(2) \rightarrow\mathscr E_0(1)\rightarrow\ 0
\end{equation}
where $\tau$ is a nowhere vanishing section of $\Omega_{\mathbb P^3}(2)$.
\end{defin}

\begin{remark}
It is possible to prove (see \cite{Barth77}, \cite{Wever81} and \cite[Example~8.4.1]{Hart2}) that $\mathscr E$ is a stable rank $2$ vector bundle on $\mathbb P^3$ with $c_1(\mathscr E)=0$ and $c_2(\mathscr E)=1$ if and only if $\mathscr E$ is isomorphic to a null-correlation bundle.
\end{remark}

\begin{ex}  \label{E}
Let $\mathscr E_0$ be a null-correlation bundle. Let $n,s\in \mathbb N$ be positive integers and let $F\colon\mathbb P^3\rightarrow \mathbb P^3$ be the absolute Frobenius on $\mathbb P^3$. Let us consider the sheaf $\mathscr E(s)={F^n}^{\star}(\mathscr E_0)\otimes \mathscr O_{\mathbb P^3}(s)$. Since $c_1({F^n}^{\star}(\mathscr E_0))=0$ and
$c_2({F^n}^{\star}(\mathscr E_0))=p^{2n}$, we see that $c_1(\mathscr E(s))=2s$ and
$c_2(\mathscr E(s))=p^{2n}+s^2$.

Let $\sigma\in H^0(\mathscr E(s))$ be a global section such that the zero locus of $\sigma$ is a curve $C$. Then we get the exact sequence:
\begin{equation} \label{S:43}
0 \rightarrow \mathscr O_{\mathbb P^3}\rightarrow \mathscr E(s) \rightarrow\mathscr I_C(2s)\rightarrow 0
\end{equation}
so that $h^0\left(\mathscr I_C(s)\right)=h^0(\mathscr E)$ and $\deg C=c_2(\mathscr E(s))=p^{2n}+s^2$. Let $H$ be a plane transversal to $C$ and $\Gamma=C\cap H$. Restricting to $H$ the exact sequence \eqref{S:43} we have:
\begin{equation}  \label{S:44}
0 \rightarrow \mathscr O_H\rightarrow \mathscr E(s)|_H \rightarrow \mathscr I_{\Gamma}(2s) \rightarrow 0
\end{equation}
so that:
\begin{equation} \label{S:45}
h^0(\mathscr I_{\Gamma}(s))=h^0(\mathscr E|_H).
\end{equation}
By \cite[Theorem 3.2]{Ein} $\mathscr E$ is stable and we can choose $H$ sufficiently general in such a way that $\mathscr E|_H$ is semi-stable, but not stable. Since $\mathscr E$ is stable and $c_1(\mathscr E)=0$, then by \cite[Lemma 3.1]{Hart1} $h^0(\mathscr E)=0$, which implies that $h^0(\mathscr I_C(s))=0$. So $C$ is not contained in any surface of degree $s$. Since $\mathscr E|_H$ is semi-stable, but not stable and $c_1(\mathscr E|_H)=0$, it must be $h^0(\mathscr E|_H)\ne 0$, so that $h^0(\mathscr I_{\Gamma}(s))\ne 0$. Moreover by \eqref{S:44} $h^0\left(\mathscr I_{\Gamma}(s-1)\right)=h^0\left(\mathscr E|_H(-1)\right)$. Now note that by \eqref{S:42} the sheaf $\mathscr E$ is determined by the exact sequence:
\begin{equation}   \label{S:46}
0\rightarrow \mathscr O_{\mathbb P^3}(-p^n)\rightarrow (F^n)^{\star}\left(\Omega_{\mathbb P^3}\right)(p^n)\rightarrow \mathscr E\rightarrow 0.
\end{equation}
so that, considered the sheaf $\mathscr F=(F^n)^{\star}(\Omega_{\mathbb P^3})$, we have the exact sequence $0\rightarrow \mathscr O_H(-p^n-1)\rightarrow\mathscr F|_H(p^n-1)\rightarrow \mathscr E|_H(-1)\rightarrow 0$, which implies that $h^0\left(\mathscr E|_H(-1)\right)=h^0\left(\mathscr F|_H(p^n-1)\right)=0$ by Lemma \ref{L:9}. So the plane curves of degree $s$ containing the generic plane section of $C$ are the minimal ones. Moreover by the previous exact sequence and by Lemma \ref{L:9} we have the equality $h^0\left(\mathscr E|_H\right)=h^0\left(\mathscr F|_H(p^n)\right)=1$, which implies by \eqref{S:45} that $h^0(\mathscr I_{\Gamma}(s))=1$. So there is a unique plane curve of degree $s$ containing $\Gamma$, which means that this plane curve of degree $s$ is the minimal plane curve containing $\Gamma$.

Now we want to know when $h^0(\mathscr E(s))\ne 0$. By \eqref{S:46} we get for each $s\in \mathbb N$:
\begin{equation} \label{S:47}
0\rightarrow \mathscr O_{\mathbb P^3}(-p^n+s)\rightarrow \mathscr F(p^n+s)\rightarrow \mathscr E(s)\rightarrow 0.
\end{equation}
By Lemma \ref{L:8} $h^0\left(\mathscr F(2p^n)\right)=6$ and $h^0\left(\mathscr F(p^n+s)\right)\ne 0$ if and only if $s\ge p^n$. So from \eqref{S:47}:
\begin{equation}  \label{S:48}
h^0\left(\mathscr E(p^n)\right)=5
\end{equation}
and $h^0\left(\mathscr E(s)\right)\ne 0$ if and only if $s\ge p^n$.

So we have global sections only for $s\ge p^n$. We want to prove there exist global sections of $\mathscr E(s)$, for every $s\ge p^n$, whose zero locus is a curve. First we must prove that the $\mathscr E$ is not split. If this was the case, then, being $\mathscr E$ a locally free sheaf of rank $2$, we would have $\mathscr E\cong \mathscr O_{\mathbb P^3}(a)\oplus\mathscr O_{\mathbb P^3}(b)$, for some $a$,~$b\in \mathbb Z$. Since $h^0(\mathscr E(p^n-1))=0$, it must be $a+p^n-1<0$ and $b+p^n-1<0$, so that $h^0(\mathscr E(p^n))\le 2$, but this contradicts \eqref{S:48}. So $\mathscr E$ is not split. Moreover, since $h^0\left(\mathscr E(p^n)\right)=5>h^0\left(\mathscr O_{\mathbb P^3}\right)=1$, by \cite[Theorem 0.1]{GRV} we get that every general nonzero global section of $\mathscr E(s)$, for $s\ge p^n$, has as zero locus a curve in $\mathbb P^3$.

Now we want to know when there are connected curves. By \cite[Proposition 1.4]{Hart2} if $h^1({\mathscr E}^{\vee}(-s))=0$, then a generic global section of $\mathscr E(s)$, for $s\ge p^r$, determines a connected curve. Note that $h^1\left({\mathscr E}^{\vee}(-s)\right)=h^2\left(\mathscr E(s-4)\right)$. By \eqref{S:47} and by Lemma \ref{L:8} $h^2\left(\mathscr E(s-4)\right)\le h^3\left(\mathscr O_{\mathbb P^3}(s-p^n-4)\right)=0$ for $s>p^n$. So for $s>p^n$ the generic global section of $\mathscr E(s)$ is connected.

Now we want to know when we have nonsingular curves. By \cite[Proposition 1.4]{Hart2}, if $\mathscr E(s-1)$ is generated by its global sections, then a sufficiently generic global section in $H^0(\mathscr E(s))$ will determine a nonsingular zero locus (not necessarily connected). Note now that in the proof of Lemma \ref{L:8} we have seen that there is a surjective morphism of sheaves ${\mathscr O}^{\oplus 6}_{\mathbb P^3}\twoheadrightarrow \mathscr F(2p^n)$ (see \eqref{S:27}). From \eqref{S:47} we see that we have also the surjective morphism $\mathscr F(2p^n)\twoheadrightarrow \mathscr E(p^n)$. So we get a surjective morphism ${\mathscr O}^{\oplus 6}_{\mathbb P^3}\twoheadrightarrow \mathscr E(p^n)$, which means that $\mathscr E(p^n)$ is generated by its global sections and so $\mathscr E(s)$ is generated by its global sections for $s\ge p^n$.

In this way we construct, for any $p$,~$n$,~$s$, with $s\ge p^n$, examples of curves $C\subset \mathbb P^3$ of degree $p^{2n}+s^2$ not contained in any surface of degree $s$ such that the minimal curve containing its generic plane section has degree $s$ and such that:
\begin{enumerate}
\item $C$ is nonsingular, in particular reduced;
\item $C$ is nonsingular and connected, which means nonsingular and irreducible, in the case $s>p^n$. In this situation the minimal curve of degree $s$ containing the generic plane section of $C$ is integral by \cite[Theorem 4.1]{Bon3}.
\end{enumerate}
In particular, we see that the bound in Theorem \ref{T:1} for connected curves is sharp. Moreover, taking $s=p^n+1$, we see that there exist connected and reduced curves (in particular nonsingular) of degree $d=2s^2-2s+1$, not lying on any surface of degree $s$, whose generic plane section is contained in an integral plane curve of degree $s$.
\end{ex}

\end{document}